\documentclass[preprint,12pt]{elsarticle}

\usepackage{amssymb}

\journal{Statistics and Probability Letters}

\newtheorem{defi}{Definition}
\newtheorem{prp}[defi]{Proposition}
\newtheorem{cor}[defi]{Corollary}

\begin{document}

\begin{frontmatter}



\title{Exponential transform of quadratic functional
and multiplicative ergodicity 
of a Gauss-Markov process\tnoteref{t1}}
\tnotetext[t1]{Dedicated to the memory of Michel Viot}



\author[MK]{Marina Kleptsyna}
\ead{Marina.Kleptsyna@univ-lemans.fr}
\address[MK]{Laboratoire Manceau de Math\'ematiques CNRS EA 3263,
Universit\'e du Maine, Avenue Olivier Messiaen, 72085 Le Mans cedex 9, France}

\author[ALB]{Alain Le Breton}
\ead{Alain.Le-Breton@imag.fr}
\author[ALB,BY2]{Bernard Ycart\corref{cor1}}
\ead{Bernard.Ycart@imag.fr}
\address[ALB]{Laboratoire Jean Kuntzmann CNRS UMR 5224,
Universit\'e Grenoble-Alpes,
51 rue des Math\'ematiques 38041 Grenoble cedex 9, France}
\address[BY2]{Laboratoire d'Excellence TOUCAN (Toulouse Cancer),
Toulouse, France}
\cortext[cor1]{corresponding author}

\begin{abstract}
The Laplace transform of partial sums of the square of
a non-centered Gauss-Markov process,
conditioning on its starting point, is explicitly computed. The parameters of
multiplicative ergodicity are deduced.
\end{abstract}


\end{frontmatter}
\section{Introduction}
Asymptotics of exponential transforms for partial sums of functionals
of a Markov chain are usually described by multiplicative ergodicity
properties, which have been thoroughly studied by Meyn and his
co-workers
\cite{BalajiMeyn00,KontoyannisMeyn03,KontoyannisMeyn05,Meyn06}; see
also \cite[p.~519]{MeynTweedie09} for a short introduction. Different
presentations of the notion, at increasingly sharp levels, can be
given. Here is the definition that will be adopted in this paper, the
notations being compatible with those of \cite{KontoyannisMeyn05}.
\begin{defi}
\label{def:multergo}
Let $\{X_t\,,\;t\in\mathbb{N}\}$ be a discrete time stochastic process,
taking values in a Polish
state space $\mathbf{X}$. Let $F$ be a measurable 
functional from the state
space $\mathbf{X}$ into $\mathbb{R}$. For $t\in\mathbb{N}$, denote by $S_t$ the
partial sum:
$$
S_t = \sum_{s=0}^{t} F(X_s)\;.
$$
The process $\{X_t\,,\;t\in \mathbb{N}\}$ is said to be
multiplicatively ergodic for 
$F$ if there exist a non-empty open subset $D$ of $\mathbb{C}$,
a function $\alpha\mapsto\Lambda(\alpha)$ from $D$ into $\mathbb{R}$, 
and a function $(\alpha,x)\mapsto
\check{f}(\alpha,x)$ from $D\times\mathbf{X}$ into $\mathbb{R}^+$, 
such that for all $\alpha\in D$ and all $x\in\mathbf{X}$: 
\begin{equation}
\label{eq:SME}
\lim_{t\to\infty}\mathbb{E}_x\left[\exp\left(\alpha S_t-t\Lambda(\alpha)
\right)\right]=\check{f}(\alpha,x)\;,
\end{equation}
where $\mathbb{E}_x$ denotes the conditional 
expectation given $X_0=x$.
\end{defi}
Balaji and Meyn \cite{BalajiMeyn00} 
introduced the notion for a Markov chain on a countable state space. It was later
extended to Markov chains on  
general state spaces by Kontoyannis and Meyn 
\cite{KontoyannisMeyn03,KontoyannisMeyn05}. Applications to large
deviations of Markov chains and Monte-Carlo simulation were developed
by Meyn \cite{Meyn06}. Since then, the notion does not seem to have
attracted much further attention. One of the reasons may be that,
except in the trivial case where the $X_t$'s are independent,
explicit calculations of the eigenvalue $\Lambda(\alpha)$ and the
eigenfunction  $\check{f}(\alpha,x)$ remain out of reach. Recently, in
the study of an exponential cell growth model
\cite{LouhichiYcart13}, the necessity of an explicit determination of
the constants of multiplicative ergodicity $\Lambda(\alpha)$ and
$\check{f}(\alpha,x)$,  was highlighted.

The main result of this note (Proposition \ref{pr:explicit})
is an explicit expression for the
exponential transform $\mathbb{E}_x[\exp(\alpha S_t)]$, in the particular
case where $\{X_t\,,\;t\in\mathbb{N}\}$ is a real-valued,
non-centered, stationary autoregressive process, and $F(x)=x^2$. The
multiplicative ergodicity coefficients $\Lambda(\alpha)$ and
$\check{f}(\alpha,x)$ are deduced (Corollary \ref{cor:explicit}). Moreover,
the convergence in (\ref{eq:SME}) is shown to be exponentially fast.
The calculation technique used for Proposition \ref{pr:explicit} was
developed in \cite{Kleptsynaetal02}, and Proposition \ref{pr:explicit}
generalizes some of the examples given in that reference. Section
\ref{statement} contains the notations and statement of the main
results, proofs are given in section \ref{proof}.
\section{Notations and statement}
\label{statement}
The process considered here is a stationary autoregressive process,
classically defined as follows.
Let $\theta$ be a real such that $0<|\theta|<1$. Let
$(e_t)_{t\geqslant 1}$ be a sequence of i.i.d.r.v.'s each following
the standard Gaussian
distribution. Let $Y_0$, independent from the sequence
$(e_t)_{t\geqslant 1}$, following the normal
$\mathcal{N}(0,(1-\theta^2)^{-1})$ distribution. For all $t\geqslant
1$ let:
$$
Y_t = \theta Y_{t-1}+e_t\;. 
$$
Then $\{Y_t\,,\;t\in\mathbb{N}\}$ 
is a stationary centered auto-regressive process.  
Denoting  by $m$ a fixed real, we consider the non-centered process
$\{X_t\,,\;t\in\mathbb{N}\}$, with $X_t=Y_t+m$. Let 
$F$ be the function $x\mapsto x^2$. 
Define:
$$
L_t(\alpha,x)=\mathbb{E}_x[\exp(\alpha S_t)] = 
\mathbb{E}_x\left[\exp\alpha\left(\sum_{s=0}^{t}F(X_s)\right)\right]\;.
$$
The particular case $m=0$ was treated in
\cite{Kleptsynaetal02}, example 4.1. The technique used here is
similar. Some notations are needed.

The explicit expression of $L_t(\alpha,x)$ uses
the two roots of the following 
equation in $\lambda$:
\begin{equation}
\label{eq:lambda}
\lambda^2-(-2\alpha+\theta^2+1)\lambda+\theta^2=0\;.
\end{equation}
Assume first that $\alpha$ is real and negative. The two roots are:
\begin{equation}
\label{lpm}
\lambda_{\pm}(\alpha) = \frac{-2\alpha+1+\theta^2\pm
\sqrt{(-2\alpha+(\theta+1)^2)(-2\alpha+(\theta-1)^2)}}{2}\;.
\end{equation}
The following inequalities hold:
\begin{equation}
\label{eq:ineq}
0<\frac{\lambda_-(\alpha)}{|\theta|}<1<\frac{\lambda_+(\alpha)}{|\theta|}\;.
\end{equation}
The two functions $\lambda_\pm(\alpha)$
admit a maximal analytic extension over an open
domain $D$ of $\mathbb{C}$, containing $(-\infty\,;\,0)$, 
over which the same inequalities hold for
their modules.
\begin{equation}
\label{eq:defD}
D = \left\{\,\alpha\in\mathbb{C}\,,\;
0<\frac{|\lambda_-(\alpha)|}{|\theta|}<1<\frac{|\lambda_+(\alpha)|}{|\theta|}
\right\}\;.
\end{equation}
In what follows, the variable is omitted in $\lambda_\pm=\lambda_\pm(\alpha)$.
Let:
\begin{equation}
\label{eq:defalphapm}
\beta_+ = \frac{1-\lambda_-}
{\lambda_+-\lambda_-}
\;;\quad
\beta_- = \frac{\lambda_+-1}
{\lambda_+-\lambda_-}\;,
\end{equation}
\begin{equation}
\label{eq:defpit}
\pi_t = \beta_+\lambda_+^{t+1}
+\beta_-\lambda_-^{t+1}\,;\;
\psi_t = \beta_+\left(\frac{\lambda_+}{\theta}\right)^{t}
+\beta_-\left(\frac{\lambda_-}{\theta}\right)^{t}\;.
\end{equation}
\begin{prp}
\label{pr:explicit}
Let:
\begin{eqnarray*}
\nu&=&\frac{m(1-\theta)}{-2\alpha+(1-\theta)^2}\,;\;
A=m(1-\theta)\nu\,;\\[2ex] 
B &=& \frac{\theta}{-2\alpha}\left(x-(1-\theta)\nu\right)^2
-\theta\nu^2
\,;\; 
C = 2\nu(x-(1-\theta)\nu)\;.
\end{eqnarray*}
Then
\begin{equation}
\label{Lumgene}
L_t(\alpha,x) = 
(\pi_t)^{-1/2}\;\exp\left(\alpha \Sigma_t\right)\;, 
\end{equation}
with
\begin{equation}
\label{sigu0}
\Sigma_t =A t+x^2+
B\left(\theta-\frac{\psi_t}{\psi_{t+1}}\right)
+C\left(\theta-\frac{1}{\psi_{t+1}}\right)
\;.
\end{equation}
\end{prp} 
Once an explicit expression of $L_t(\alpha,x)$ has been obtained,
deriving its asymptotics as $t$ tends to infinity is easy, using 
(\ref{eq:ineq}). From the same inequalities, it follows that
the convergence in (\ref{eq:SME}) holds at exponential speed
$O((\theta/\lambda_+)^t)$. 
\begin{cor}
\label{cor:explicit}
For all $m\in \mathbb{R}$,
the process $\{X_t\,,\;t\in\mathbb{N}\}$ is multiplicatively ergodic for $F$
in the sense of Definition \ref{def:multergo}. 
The domain $D$ is defined by (\ref{eq:defD}). For $\alpha\in D$ and
$x\in\mathbb{R}$, the limit
(\ref{eq:SME}) holds with
\begin{equation}
\label{defL}
\Lambda(\alpha)=\frac{\alpha
    m^2(1-\theta)^2}{(-2\alpha+(1-\theta)^2)} - \frac{1}{2}\log(\lambda_+)
\;,
\end{equation}
and 
\begin{equation}
\label{defcheckf}
\check{f}(\alpha,x) =
(\beta_+\lambda_+)^{-1/2}
\exp\left(\alpha\left(x^2+B\left(\theta-\frac{\theta}{\lambda_+}
\right)+C\theta\right)\right)\;. 
\end{equation}
\end{cor}
In (\ref{defL}), $\log$ denotes an analytic extension to $D$ of the
ordinary logarithm on $\mathbb{R}^+$ (recall that $\lambda_+$ is positive for 
$\alpha\in (-\infty\,;0\,)\subset D$).
\section{Proof}
\label{proof}
The technique of proof is an application of Theorem 1 in
\cite{Kleptsynaetal02}. As in that reference,
we shall rewrite the exponential transform 
$\mathbb{E}_x[\exp(\alpha S_t)]$ as a Laplace transform, and use $\mu/2$ as
its variable: in what follows, $\alpha$ is replaced by $\mu =
-2\alpha$, and $L_t(\alpha,x)$ by $L_t(\mu,x)$. 
Since by definition $X_t=Y_t+m$, the initial condition $X_0=x$ amounts to
$Y_0=x-m$. Hence the conditional distribution of $\{X_t\,,\,t\in\mathbb{N}\}$
given $X_0=x$ coincides with the distribution of
$\{Y^x_t+m\,,\;t\in\mathbb{N}\}$,
where $\{Y^x_t\,,\;t\in\mathbb{N}\}$ is
defined by $Y^x_0=x-m$, and for all $t\geqslant 1$,
$$
Y^x_t = \theta Y^x_{t-1} +e_t\;.
$$
So our aim is to compute
$$
L_t(\mu,x)
=
\mathbb{E}
\left[\exp\left(-\frac{\mu}{2}\sum_{s=0}^{t} (Y_s^x+m)^2\right)\right]\;.
$$ 
The process $\{Y^x_t+m\,,\;t\geqslant 0\}$ has the same covariance function
as $\{Y_t^x\,,\;t\geqslant 0\}$. 
Denote by $m^x_t$ its mean function. It is such that 
$m^x_0=x$ and for $s\geqslant 1$,
\begin{equation}
\label{recm}
m_s^x=\theta m_{s-1}^x+m(1-\theta)\;.
\end{equation}
Theorem 1 of \cite{Kleptsynaetal02}, and the
calculations in
Example 4.1 therein, yield
(\ref{Lumgene}), where $\pi_t$ is defined by (\ref{eq:defpit}), and
\begin{equation}
\label{sigu1}
\Sigma_t = \sum_{s=0}^t \frac{\psi_s}{\theta\psi_{s+1}} (z_s)^2\;,
\end{equation}
where $z_0=m_0^x=x$ and
\begin{equation}
\label{zeps}
z_s = m^x_s-\sum_{r=0}^{s-1} \theta^{s-r}
\left(1-\frac{\psi_r}{\theta\psi_{r+1}}\right) z_r\;.
\end{equation}
Using the expression of $m^x_s$ (\ref{recm}):
\begin{eqnarray*}
z_s &=& (1\!-\!\theta)m+\theta m^x_{s-1}
-\theta\! \left(1-\frac{\psi_{s-1}}{\theta\psi_{s}}\right) z_{s-1} 
-\theta\!\sum_{r=0}^{s-2} \theta^{s-1-r}
\left(1-\frac{\psi_r}{\theta\psi_{r+1}}\right) z_r\\[2ex]
&=&(1-\theta)m +\theta z_{s-1}-\theta 
\left(1-\frac{\psi_{s-1}}{\theta\psi_{s}}\right) z_{s-1}
=\frac{\psi_{s-1}}{\psi_{s}} z_{s-1}+ (1-\theta)m\;.
\end{eqnarray*}
The resolvent of that equation is $(\psi_s)^{-1}$, from which the following
expression is obtained.
$$
z_s = \psi_s^{-1}
\left(x+m(1-\theta)\sum_{l=1}^s\psi_l\right)\;. 
$$
Plugging into (\ref{sigu1}) yields:
\begin{equation}
\label{sigu2}
\Sigma_t = \sum_{s=0}^t \frac{1}{\theta \psi_s \psi_{s+1}}
\left(x+m(1-\theta)\sum_{l=1}^s\psi_l\right)^2\;.
\end{equation}
Setting $z=\lambda_+/\theta$, $z^{-1}=\lambda_-/\theta$, one has
$\psi_s = \beta_+z^s+\beta_-z^{-s}$, and
$$
\sum_{l=1}^s\psi_l= \beta_+\frac{z}{1-z}(1-z^s)
+\beta_-\frac{z^{-1}}{1-z^{-1}}(1-z^{-s})\;.
$$
Define
$$
\Delta_s = x+m(1-\theta)\sum_{l=1}^s\psi_l 
= a_+z^s +a_-z^{-s}+a \;.
$$
with
$$
a_+=m(\theta-1)\beta_+\frac{z}{1-z}\;;\quad
a_-=m(\theta-1)\beta_-\frac{z^{-1}}{1-z^{-1}}\;,
$$
and $a=x-(a_++a_-)$.
For $s\geqslant 1$,
$$
\psi_{s+1}\psi_{s-1}-(\psi_{s})^2
=
\beta_+\beta_-(z-z^{-1})^2\;,
$$
hence
$$
\frac{\psi_{s-1}}{\psi_{s}}
-\frac{\psi_{s}}{\psi_{s+1}}
=\frac{\beta_+\beta_-(z-z^{-1})^2}
{\psi_s\psi_{s+1}}\;.
$$
Now let us choose the three constants $A$, 
$B$, $C$ such that:
\begin{equation}
\label{deltasu}
\frac{(\Delta_s)^2}{\theta\psi_s\psi_{s+1}}
=
A+B\left(\frac{\psi_{s-1}}{\psi_{s}}
-\frac{\psi_{s}}{\psi_{s+1}}
\right)
+C\left(\frac{1}{\psi_{s}}
-\frac{1}{\psi_{s+1}}
\right)\;.
\end{equation}
For that we need:
$$
(\Delta_s)^2=A\theta\psi_s\psi_{s+1}
+B\theta\beta_+\beta_-(z-z^{-1})^2+
C\theta(\psi_{s+1}-\psi_s)\;,
$$
then
\begin{eqnarray*}
&&(a_+z^s +a_-z^{-s}+a)^2=
A\theta(\beta_+z^s +\beta_-z^{-s})(\beta_+z^{s+1}+\beta_-z^{-s-1})\\
&&\hspace*{2cm}+B\theta\beta_+\beta_-(z-z^{-1})^2
+C\theta(\beta_+(z-1)z^s+\beta_-(z^{-1}-1)z^{-s})\;.
\end{eqnarray*}
The expressions of $A$, $B$, $C$ can be guessed by identifying
powers of $z$ in the expression above.
The following constants satisfy the requirements.
\begin{eqnarray*}
A&=&\displaystyle{\frac{(a_+)^2}{\theta(\beta_+)^2z}\;,}\\[2ex]
B&=&\displaystyle{\frac{a^2+2a_+a_--A\theta\beta_+\beta_-(z+z^{-1})}
{\theta\beta_+\beta_-(z-z^{-1})^2}\;,}\\[2ex]
C&=&\displaystyle{\frac{2aa_+}{2\beta_+(z-1)}\;.}
\end{eqnarray*}
For these constants, plugging (\ref{deltasu}) into (\ref{sigu2})
and summing yields:
\begin{equation}
\label{sigu3}
\Sigma_t
= \frac{x^2}{\theta\psi_0\psi_1} + At
+ B\left(\frac{\psi_0}{\psi_1}-\frac{\psi_t}{\psi_{t+1}}\right)
+ C\left(\frac{1}{\psi_1}-\frac{1}{\psi_{t+1}}\right)\;.
\end{equation}
Substituting $\psi_0=1$ and $\psi_1=\theta^{-1}$
in (\ref{sigu3}) gives (\ref{sigu0}).
\vskip 2mm\noindent
To finish the proof of Proposition \ref{pr:explicit}, more explicit
expressions of $A$, $B$, and $C$ must be obtained.
As a preliminary observation, recall that $z=\lambda_+/\theta$
and $z^{-1}=\lambda_-/\theta$ are the two roots of
$
z^2-\theta^{-1}(\mu+\theta^2+1)z+1=0\;.
$
From there the following symmetric functions of the two roots are
obtained.
\begin{eqnarray*}
z+z^{-1}&=&\theta^{-1}(\mu+\theta^2+1)\;,\\
(z-z^{-1})^2&=&\theta^{-2}(\mu+(1-\theta)^2)(\mu+(1+\theta)^2)\;,\\
(z-1)(1-z^{-1})&=&z+z^{-1}-2=\theta^{-1}(\mu+(1-\theta)^2)\;,\\
(z+1)(1+z^{-1})&=&z+z^{-1}+2=\theta^{-1}(\mu+(1+\theta)^2)\;.
\end{eqnarray*}
The constants $\beta_+$ and $\beta_-$ can be written as 
functions of $z$:
$$
\beta_+ = \frac{1-\theta z^{-1}}{\theta(z-z^{-1})}
\quad\mbox{and}\quad
\beta_- = \frac{\theta z -1}{\theta(z-z^{-1})}\;.
$$
From there:
\begin{eqnarray*}
\beta_+\beta_-
&=&\displaystyle{\frac{(1-\theta z^{-1})(\theta z -1)}{\theta^2(z-z^{-1})^2}}
=\displaystyle{\frac{\theta(z+z^{-1})-1-\theta^2}{\theta^2(z-z^{-1})^2}}\\[2ex]
&=&\displaystyle{\frac{\mu+\theta^2+1-1-\theta^2}{\theta^2(z-z^{-1})^2}}
=\displaystyle{\frac{\mu}{(\mu+(1-\theta)^2)(\mu+(1-\theta)^2)}\;.}
\end{eqnarray*}
Now:
\begin{eqnarray*}
A&=&\displaystyle{\frac{a_+^2}{\theta\beta_+^2 z}}
=\displaystyle{\frac{m^2(1-\theta)^2z(1-z^{-1})^2}
{\theta((z-1)(1-z^{-1}))^2}}\\[2ex]
&=&\displaystyle{\frac{m^2(1-\theta)^2(z+z^{-1}-2)}
{\theta((z-1)(1-z^{-1}))^2}}
=\displaystyle{\frac{m^2(1-\theta)^2}{\mu+(1-\theta)^2}
= m(1-\theta)\nu\;.}
\end{eqnarray*}
Here is the calculation of $B$:
$$
B=\frac{a^2+2a_+a_--A\theta\beta_+\beta_-(z+z^{-1})}
{\theta\beta_+\beta_-(z-z^{-1})^2}=\frac{N}{D}\;.
$$
From the expression of $\beta_+\beta_-$ above, $D=\mu/\theta$. The numerator 
$N$ contains three terms. The first term is:
$$
a^2=(x-(a_++a_-))^2\;,
$$
where,
\begin{eqnarray*}
a_++a_-
&=& \displaystyle{-m(1-\theta)\left(
\beta_+\frac{z}{1-z}+\beta_-\frac{z^{-1}}{1-z^{-1}}\right)\;,}\\[2ex]
&=& \displaystyle{-m(1-\theta)\frac{
\beta_+z(1-z^{-1})+\beta_- z^{-1}(1-z)}{(1-z)(1-z^{-1})}\;,}\\[2ex]
&=& \displaystyle{-m(1-\theta)
\frac{\psi_1-\psi_0}{(1-z)(1-z^{-1})}\;,}\\[2ex]
&=& \displaystyle{\frac{m(1-\theta)^2}
{\mu+(1-\theta)^2}=(1-\theta)\nu\;.}
\end{eqnarray*}
Thus:
$$
a^2=(x-(1-\theta)\nu)^2\;.
$$
The second term in the numerator $N$ is:
\begin{eqnarray*}
2a_+a_-&=&\displaystyle{
\frac{2m^2(1-\theta)^2\beta_+\beta_-}{(1-z)(1-z^{-1})}
}\\[2ex]
&=&\displaystyle{
\frac{2m^2(1-\theta)^2(-\theta\mu)}{(\mu+(1+\theta)^2)(\mu+(1-\theta)^2)^2}
}\\[2ex]
&=&\displaystyle{-\frac{2\theta\mu}{\mu+(1+\theta)^2}}\,\nu^2\;.
\end{eqnarray*}
The last term in the numerator $N$ is:
\begin{eqnarray*}
-\theta A\beta_+\beta_-(z+z^{-1})
&=&\displaystyle{-\theta A
\frac{\mu(\mu+\theta^2+1)}{(\mu+(1+\theta)^2)(\mu+(1-\theta)^2)}
}\\[2ex]
&=&\displaystyle{-
\frac{m^2(1-\theta)^2\mu(\mu+\theta^2+1)}{(\mu+(1+\theta)^2)(\mu+(1-\theta)^2)^2}
}\\[2ex]
&=&\displaystyle{
-\frac{\mu(\mu+\theta^2+1)}{\mu+(1+\theta)^2}\,\nu^2
\;.}
\end{eqnarray*}
Grouping the three terms and multiplying by $D^{-1}=\theta/\mu$:
\begin{eqnarray*}
B&=&\displaystyle{\frac{\theta}{\mu}(x-(1-\theta)\nu)^2
-\frac{2\theta^2+\theta(\mu+\theta^2+1)}{\mu+(1+\theta)^2}\,\nu^2\;,}\\[2ex]
&=&\displaystyle{\frac{\theta}{\mu}(x-(1-\theta) \nu)^2
-\theta\nu^2\;.}
\end{eqnarray*}
Finally, here is the calculation of $C$.
\begin{eqnarray*}
C&=&\displaystyle{\frac{2aa_+}{\theta\beta_+(z-1)}}\\[2ex]
&=&\displaystyle{2(x-(1-\theta)\nu)\times
\frac{m(1-\theta)z}{\theta(z-1)^2}\;,}
\end{eqnarray*}
where
\begin{eqnarray*}
\displaystyle{\frac{z}{\theta(z-1)^2}}
&=&\displaystyle{\frac{z(1-z^{-1})^2}{\theta((z-1)(1-z^{-1}))^2}
}\\[2ex]
&=&\displaystyle{\frac{z+z^{-1}-2}{\theta((z-1)(1-z^{-1}))^2}}\\[2ex]
&=&\displaystyle{\frac{1}{\mu+(1-\theta)^2}\;.}
\end{eqnarray*}
Therefore:
$$
C=2\nu(x-(1-\theta)\nu)\;.
$$
This ends the proof of Proposition \ref{pr:explicit}.

\noindent\textbf{Concluding remarks}
\begin{itemize}
\item[(a)] The (non-conditional) exponential transform
  $\mathbb{E}[\exp(\alpha S_t)]$ can be obtained by integrating the
  right-hand side of (\ref{Lumgene}) against the centered Gaussian
  distribution with variance $(1-\theta^2)^{-1}$. A
  simpler calculation can be carried through, applying again Theorem 1 of
  \cite{Kleptsynaetal02} as was done above. Similar explicit
  expressions are obtained, that will not be reproduced here.
\item[(b)]
Here, the parameters of multiplicative ergodicity have been deduced
from the explicit expression of the Laplace transform,
for a specific class of Gaussian processes. On the basis of
the results of \cite{Kleptsynaetal02}, there is some hope that they
could be directly calculated for more general processes. This should
be investigated in a forthcoming paper.
\end{itemize}








\end{document}